\documentclass[11pt]{article}
\usepackage{graphics}
\usepackage{amsfonts}
\usepackage{amsmath}
\usepackage{amsthm}
\usepackage{amssymb}
\usepackage{graphicx}

\newtheorem{theorem}{Theorem}[section]
\newtheorem{lemma}{Lemma}[section]
\newtheorem{remark}{Remark}[section]

\input{tcilatex}

\begin{document}

\title{BDDC and FETI-DP under Minimalist Assumptions}
\author{Jan Mandel{}\footnotemark[1]  \footnotemark[3]  \and Bed\v{r}ich
Soused\'{\i}k{}\footnotemark[1]  \footnotemark[2]  \footnotemark[3]
\footnotemark[4]  }
\maketitle

\begin{abstract}
The FETI-DP, BDDC and P-FETI-DP preconditioners are derived in a particulary
simple abstract form. It is shown that their properties can be obtained from
only on a very small set of algebraic assumptions. The presentation is
purely algebraic and it does not use any particular definition of method
components, such as substructures and coarse degrees of freedom. It is then
shown that P-FETI-DP and BDDC are in fact the same. The FETI-DP and the BDDC
preconditioned operators are of the same algebraic form, and the standard
condition number bound carries over to arbitrary abstract operators of this
form. The equality of eigenvalues of BDDC and FETI-DP also holds in the
minimalist abstract setting. The abstract framework is explained on a
standard substructuring example.
\end{abstract}

\footnotetext[1]{%
Department of Mathematical Sciences, University of Colorado at Denver and
Health Sciences Center, P.O. Box 173364, Campus Box 170, Denver, CO 80217,
USA}

\footnotetext[2]{%
Department of Mathematics, Faculty of Civil Engineering, Czech Technical
University, Th\'akurova 7, 166 29 Prague 6, Czech Republic}

\footnotetext[3]{%
Supported by the National Science Foundation under grants CNS-0325314,
CNS-0719641, and DMS-0713876.}

\footnotetext[4]{%
Supported by the program of the Information society of the Academy of
Sciences of the Czech Republic \mbox{1ET400760509} and by the research
project \mbox{CEZ MSM 6840770003}.}

\section{Introduction}

\label{sec:introduction}

The BDDC and FETI-DP methods are iterative substructuring methods that use
coarse degrees of freedom associated with corners, edges or faces between
subdomains, and they are currently the most advanced versions of the BDD and
FETI families of methods. The BDDC method introduced by Dohrmann~\cite%
{Dohrmann-2003-PSC} is a Neumann-Neumann method of Schwarz type~\cite%
{Dryja-1995-SMN}. The BDDC method iterates on the system of primal
variables\ reduced to the interfaces between the subdomains, and it can be
understood as further development of the BDD method by Mandel~\cite%
{Mandel-1993-BDD}. The FETI-DP method by Farhat et al.~\cite%
{Farhat-2001-FDP,Farhat-2000-SDP} is a dual method that iterates on a system
for Lagrange multipliers that enforce continuity on the interfaces, and it
is a further development of the FETI method by Farhat and Roux~\cite%
{Farhat-1991-MFE}. Algebraic relations between FETI and BDD methods were
pointed out by Rixen et al.~\cite{Rixen-1999-TCF}, Klawonn and Widlund~\cite%
{Klawonn-2001-FNN}, and Fragakis and Papadrakakis~\cite{Fragakis-2003-MHP}.
A common bound on the condition number of both the FETI and the BDD method
in terms of a single inequality in was given in~\cite{Klawonn-2001-FNN}. In
the case of corner constraints only, a method same as BDDC was suggested by
Cros~\cite{Cros-2003-PSC}. Fragakis and Papadrakakis~\cite{Fragakis-2003-MHP}
derived primal versions of FETI and FETI-DP, called respectively P-FETI and P-FETI-DP, and they have also observed
that the eigenvalues of BDD and a certain version of FETI are identical.
Mandel, Dohrmann, and Tezaur~\cite{Mandel-2005-ATP} have proved that the
eigenvalues of BDDC and FETI-DP are identical and they have obtained a
simplifed and fully algebraic version (i.e., with no undetermined constants)
of a common condition number estimate for BDDC and FETI-DP, similar to the
estimate by Klawonn and Widlund~\cite{Klawonn-2001-FNN} for BDD and FETI.
Simpler proofs of the equality of eigenvalues of BDDC and FETI-DP were
obtained by Li and Widlund~\cite{Li-2006-FBB}, and by Brenner and Sung~\cite%
{Brenner-2007-BFW}, who also gave an example when BDDC\ has an eigenvalue
equal to one but FETI-DP\ does not. The proof of the equality of eigenvalues
of BDD and a certain version of FETI and of FETI-DP and P-FETI-DP was recently given by Fragakis~%
\cite{Fragakis-2007-FDD} in a more general framework.

In this contribution, we derive the FETI-DP, BDDC and P-FETI-DP\
preconditioners in a particulary simple abstract form, with only a very
small set of algebraic assumptions (Sec. \ref{sec:formulation}). The
presentation is purely algebraical and it does not use any particular
definition of method components, such as substructures and coarse degrees of
freedom. We show that P-FETI-DP\ and BDDC are in fact the same. We then
present the condition number bound and the proof of the equality of
eigenvalues of BDDC and FETI-DP, in the minimalist abstract setting (Sec. %
\ref{sec:connections}). Finally, we illustrate the abstract framework on a
substructuring example (Sec. \ref{sec:substructuring}).

\section{Notation and Preliminaries}

\label{sec:preliminaries}All spaces in this paper are finite dimensional
linear spaces. The dual space of a space $V$\ is denoted by $V^{\prime}$ and
$\left\langle \cdot,\cdot\right\rangle $\ is the duality pairing. For a
linear operator $L:W\rightarrow V$ we define its transpose $L^{T}:V^{\prime
}\rightarrow W^{\prime}$ by $\left\langle v,Lw\right\rangle =\left\langle
L^{T}v,w\right\rangle $ for all $v\in V^{\prime}$, $w\in W$,\ and $%
\left\Vert v\right\Vert _{K}=\sqrt{\left\langle Kv,v\right\rangle }$ denotes
the norm associated with a symmetric and positive definite operator $%
K:V\rightarrow V^{\prime}$, i.e., such that $\left\langle Kv,v\right\rangle
>0$ for all $v\in V$, $v\neq0$. The norm of a linear operator $%
E:V\rightarrow V$ subordinate to this vector norm is defined by $\left\Vert
E\right\Vert _{K}=\max_{v\in V,v\neq0}\left\Vert Ev\right\Vert
_{K}/\left\Vert v\right\Vert _{K}$. The notation $I_{V}$ denotes the
identity operator on the space $V$.

Mappings from a space to its dual arise naturally in the variational setting
of systems of linear algebraic equations. An an example, consider an $%
n\times n$ matrix $A$ and the system of equations $Ax=b$. The variational
form of this system is%
\begin{equation*}
x\in V:\left( Ax,y\right) =\left( b,y\right) \quad\forall y\in V,
\end{equation*}
where $V=\mathbb{R}^{n}$ and $\left( \cdot,\cdot\right) $ is the usual
Euclidean inner product on $\mathbb{R}^{n}$. For a fixed $x$, instead of the
value $Ax$, we find it convenient to consider the linear mapping $%
y\mapsto\left( Ax,y\right) $. This mapping is an element of the dual space $%
V^{\prime}$. Denote this mapping by $Kx$ and its value at $y$ by $%
\left\langle Kx,y\right\rangle $; then $K:V\rightarrow V^{\prime}$ is a
linear operator from $V$ to its dual that corresponds to $A$. This setting
involving dual spaces is convenient and compact when dealing with multiple
nested spaces, or with dual methods (such as FETI). Restricting a linear
functional to a subspace is immediate, while the equivalent notation without
duality requires introducing new operators, namely projections or transposes
of injections. Also, this setting allows us to make a clear distinction
between an approximate solution and its residual, which is in the dual
space. It is beneficial to have approximate solutions and residuals in
different spaces, because they need to be treated differently.

We wish to solve a system of linear algebraic equations%
\begin{equation*}
Ku=f,
\end{equation*}
where $K:V\rightarrow V^{\prime}$, by a preconditioned conjugate gradient
method. Here, a preconditioner is a mapping $M:V^{\prime}\rightarrow V$. In
iteration $k$ the method computes the residual%
\begin{equation*}
r^{(k)}=Ku^{(k)}-f\in V^{\prime},
\end{equation*}
and the preconditioner computes the increment\ to the approximate solution $%
u^{(k)}$ as a linear combination of the preconditioned residual $Mr^{(k)}\in
V$ with preconditioned residuals in earlier iterations. Convergence
properties of the method can be established from the eigenvalues $\lambda$\
of the \emph{preconditioned operator} $MK$; the condition number%
\begin{equation*}
\kappa=\frac{\lambda_{\max}(MK)}{\lambda_{\min}(MK)},
\end{equation*}
gives a well-known bound on the error reduction, cf. e.g.~\cite%
{Golub:1989:MAC},
\begin{equation*}
\bigl\Vert e^{(k)}\bigr\Vert _{K}\leq2\left( \frac{\sqrt{\kappa}-1}{\sqrt{%
\kappa}+1}\right) ^{k}\bigl\Vert e^{(0)}\bigr\Vert _{K},
\end{equation*}
where $e^{(k)}=u^{(k)}-u$ is the error of the solution in iteration $k$.

\section{Abstract Formulation of the Preconditioners}

\label{sec:formulation}

\subsection{Setting and Assumptions}

We now list a minimalist set of spaces, linear operators, and assumptions
needed to formulate the BDDC and the FETI-DP methods and to prove their
properties. Let $W$ be a finite dimensional space and let $a\left(
\cdot,\cdot\right) $ be a symmetric positive semi-definite\ bilinear form on
$W$. Let $\widehat{W}\subset W$ be a subspace such that $a$ is positive
definite on $\widehat{W}$, and $f\in\widehat{W}^{\prime}$. We wish to solve
a variational problem%
\begin{equation}
u\in\widehat{W}:a\left( u,v\right) =\left\langle f,v\right\rangle
\quad\forall v\in\widehat{W}.  \label{eq:variational}
\end{equation}

The preconditioners we are interested in are characterized by a selection of
an intermediate space $\widetilde{W}$,%
\begin{equation}
\widehat{W}\subset\widetilde{W}\subset W,  \label{eq:w-tilde}
\end{equation}
such that
\begin{equation}
a\left( \cdot,\cdot\right) \text{ is positive definite on }\widetilde{W},
\label{eq:pos-def}
\end{equation}
and a selection of linear operators $E$, $B$, and $B_{D}$. The operator $E$
is a projection onto $\widehat{W}$,
\begin{equation}
E:\widetilde{W}\rightarrow\widehat{W},\quad E^{2}=E,\quad\limfunc{range}E=%
\widehat{W}.  \label{eq:E}
\end{equation}
The role of the operator $B$ is to enforce the condition $u\in\widehat{W}$
by $Bu=0$,
\begin{equation}
B:\widetilde{W}\rightarrow\Lambda,\quad\limfunc{null}B=\widehat{W},\quad%
\limfunc{range}B=\Lambda.  \label{eq:BL}
\end{equation}
The operator $B_{D}^{T}$ is a generalised inverse of $B$,
\begin{equation}
B_{D}^{T}:\Lambda\rightarrow\widetilde{W},\quad BB_{D}^{T}=I_{\Lambda}.
\label{eq:BD}
\end{equation}
\emph{The properties (\ref{eq:w-tilde}) -- (\ref{eq:BD}) are enough for the
BDDC and the FETI-DP theories separately, and they will be assumed from now
on.} To relate the two methods, we shall also assume that
\begin{equation}
B_{D}^{T}B+E=I  \label{eq:corresp}
\end{equation}
\emph{when needed.}\/ No further assumptions are made in the rest of the
paper.

To formulate the preconditioners and their properties, we need to define
several more linear operators from the concepts already introduced. Denote
by
\begin{equation}
R:\widehat{W}\rightarrow\widetilde{W},\quad R:w\in\widehat{W}\longmapsto w\in%
\widetilde{W},  \label{eq:def-r}
\end{equation}
the natural injection from $\widehat{W}$ to $\widetilde{W}$. Clearly,%
\begin{equation}
ER=I_{\widehat{W}}.  \label{eq:ER}
\end{equation}

\begin{remark}
The conditions (\ref{eq:pos-def}) -- (\ref{eq:corresp}) are satisfied in the
applications of FETI-DP and BDDC methods, see \cite{Mandel-2005-ATP} and
Sec. \ref{sec:substructuring}. Note that the assumption (\ref{eq:BD}) allows
the case when $B$ is a matrix that does not have full row rank. All that is
needed is to define $\Lambda$ as $\limfunc{range}B$. In the literature, \cite%
{Mandel-2005-ATP} and references therein, the projection $E$ is often
written in the form $E=RD_{P}R^{T}$ where $R$ is a mapping of another space
(isomorphic to $\widehat{W}$) into $W$. In the abstract setting here, we
choose to formulate the methods in the space $\widehat{W}$ directly, it
turns out that the space $W$ is not needed for the theory at all, and $R$
becomes the identity embedding of $\widehat{W}$ into $\widetilde{W}$. The
equation (\ref{eq:BD}) is found already in \cite[Lemma 1]{Rixen-1999-TCF} in
a special case. It was extended to form used presently and to cover more
general algorithms, and used to obtain important connections between dual
and primal substructuring methods in \cite%
{Fragakis-2003-MHP,Klawonn-2001-FNN}.
\end{remark}

We also define the linear operators $\widehat{S}$ and $\widetilde{S}$
associated with the bilinear form $a$ on the spaces $\widehat{W}$ and $%
\widetilde{W}$, respectively, by%
\begin{align}
\widehat{S} & :\widehat{W}\rightarrow\widehat{W}^{\prime},\quad\left\langle
\widehat{S}v,w\right\rangle =a\left( v,w\right) \quad\forall v,w\in \widehat{%
W}.  \label{eq:def-s-hat} \\
\widetilde{S} & :\widetilde{W}\rightarrow\widetilde{W}^{\prime},\quad\left%
\langle \widetilde{S}v,w\right\rangle =a\left( v,w\right) \quad\forall v,w\in%
\widetilde{W},  \label{eq:def-s-tilde}
\end{align}
From (\ref{eq:def-s-hat}), the variational problem~(\ref{eq:variational})
becomes%
\begin{equation}
\widehat{S}u=f.  \label{eq:op-form}
\end{equation}
Further, it follows from (\ref{eq:def-s-hat}), (\ref{eq:def-s-tilde}), and (%
\ref{eq:def-r}), that%
\begin{equation}
\widehat{S}=R^{T}\widetilde{S}R.  \label{eq:assembly-hat}
\end{equation}

\subsection{BDDC}

The following presentation of BDDC\ follows \cite{Mandel-2007-ASF}. It is
essentially same as the approach of \cite{Brenner-2007-BFW}, and related to
the concept of subassembly in~\cite{Li-2006-FBB}. The \emph{BDDC} \cite%
{Dohrmann-2003-PSC} is the method of preconditioned conjugate gradients
applied to the system (\ref{eq:op-form}), with the abstract BDDC
preconditioner $M_{BDDC}:\widehat{W}^{\prime}\rightarrow\widehat{W}$ defined
as%
\begin{equation}
M_{BDDC}:r\longmapsto u=Ew,\quad w\in\widetilde{W}:\quad a\left( w,z\right)
=\left\langle r,Ez\right\rangle ,\quad\forall z\in\widetilde{W}.
\label{eq:variational-BDDC}
\end{equation}
For the equivalence of (\ref{eq:variational-BDDC}) with other formulations
of BDDC, see \cite[Lemma 7]{Mandel-2005-ATP}.

From the definitions of $\widetilde{S}$ in (\ref{eq:def-s-tilde}) and $R$ in
(\ref{eq:def-r}), it follows that the operator form of the BDDC
preconditioner is%
\begin{equation}
M_{BDDC}=E\widetilde{S}^{-1}E^{T}.  \label{eq:operator-bddc}
\end{equation}

\subsection{FETI-DP}

This presentation of FETI-DP\ follows \cite{Mandel-2005-ATP}. The
variational problem~(\ref{eq:variational}) is equivalent to the minimization
\begin{equation}
\frac{1}{2}a\left( u,u\right) -\left\langle f,u\right\rangle \rightarrow \min%
\text{ subject to }u\in\widehat{W}.  \label{eq:minimization}
\end{equation}
Using
\begin{equation*}
\left\langle f,u\right\rangle =\left\langle f,Eu\right\rangle =\left\langle
E^{T}f,u\right\rangle ,\quad u\in\widehat{W},
\end{equation*}
we can write (\ref{eq:minimization}) as a constrained minimization problem
posed on $\widetilde{W}$,%
\begin{equation*}
\frac{1}{2}a\left( u,u\right) -\left\langle E^{T}f,u\right\rangle
\rightarrow\min\text{ subject to }u\in\widetilde{W}\quad\text{and}\quad Bu=0.
\end{equation*}
Introducing the Lagrangean
\begin{equation*}
{\mathcal{L}(w,\lambda)=}\frac{1}{2}a\left( u,u\right) -\left\langle
E^{T}f,u\right\rangle +\left\langle B^{T}\lambda,u\right\rangle ,
\end{equation*}
where $\lambda\in\Lambda^{\prime}$ are the Lagrange multipliers, we obtain
that problem (\ref{eq:minimization}) is equivalent to solving the
saddle-point problem \cite{Mangasarian-1994-NP}
\begin{equation}
{\min\limits_{w\in\widetilde{W}}\max\limits_{\lambda\in\Lambda^{\prime}}%
\mathcal{L}(w,\lambda)}.  \label{eq:minimax}
\end{equation}
Since
\begin{equation*}
{\min\limits_{w\in\widetilde{W}}\max\limits_{\lambda\in\Lambda^{\prime}}%
\mathcal{L}(w,\lambda)=\max\limits_{\lambda\in\Lambda^{\prime}}\min%
\limits_{w\in\widetilde{W}}\mathcal{L}(w,\lambda),}
\end{equation*}
it follows that (\ref{eq:minimization}) is equivalent to the dual problem
\begin{equation}
\frac{\partial\mathcal{F}(\lambda)}{\partial\lambda}=0,
\label{eq:dual-problem}
\end{equation}
where
\begin{equation}
\mathcal{F}(\lambda)=\min\limits_{w\in\widetilde{W}}\mathcal{L}(w,\lambda).
\label{eq:dual-functional}
\end{equation}

Problem (\ref{eq:dual-problem}) is equivalent to stationary conditions for
the Lagrangean $\mathcal{L}$,%
\begin{align}
\frac{\partial}{\partial w}\mathcal{L}(w,\lambda) & \perp\widetilde {W},
\label{eq:stationary-w} \\
\frac{\partial}{\partial\lambda}\mathcal{L}(w,\lambda) & =0,  \notag
\end{align}
which is the same as solving for $w\in\widetilde{W}$ and $\lambda\in
\Lambda^{\prime}$ from the system
\begin{equation}
\begin{array}{ccccc}
\widetilde{S}w & + & B^{T}\lambda & = & E^{T}f, \\
Bw &  &  & = & 0.%
\end{array}
\label{eq:FETI-DP-system}
\end{equation}
Solving from the first equation in (\ref{eq:FETI-DP-system}) and
substituting into the second equation, we get the dual problem in an
operator form,%
\begin{equation}
B\widetilde{S}^{-1}B^{T}\lambda=B\widetilde{S}^{-1}E^{T}f.
\label{eq:FETI-DP-system-2}
\end{equation}

The \emph{FETI-DP method} \cite{Farhat-2001-FDP,Farhat-2000-SDP} is the
method of preconditioned conjugate gradients applied to the problem (\ref%
{eq:FETI-DP-system-2}), with the preconditioner given by
\begin{equation}
M_{FETI-DP}=B_{D}\widetilde{S}B_{D}^{T}.  \label{eq:Dirichlet-preconditioner}
\end{equation}

The FETI-DP method solves for the Lagrange multiplier $\lambda$. The
corresponding primal solution is found as the minimizer of $w$ in (\ref%
{eq:dual-functional}), or equivalently, from~(\ref{eq:stationary-w}), which
is the same as the first equation in (\ref{eq:FETI-DP-system}); hence,%
\begin{equation*}
w=\widetilde{S}^{-1}\left( E^{T}f-B^{T}\lambda\right) .
\end{equation*}
If $\lambda$ is the exact solution of the dual problem (\ref{eq:dual-problem}%
), then $w\in\widehat{W}$ and so $u=w$ is the desired solution of the primal
minimization problem (\ref{eq:minimization}). However, for approximate
solution $\lambda$, in general $w\notin\widehat{W}$, and so the primal
solution needs to be projected onto $\widehat{W}.$ We use the operator $E$
for this purpose. So, for an arbitrary Lagrange multiplier $\lambda$, the
corresponding approximate solution of the original problem is
\begin{equation}
u=E\widetilde{S}^{-1}\left( E^{T}f-B^{T}\lambda\right) .
\label{eq:primal-sol}
\end{equation}
Note that the operator $E$ does not play any role in FETI-DP iterations
themselves. It only serves to form the right-hand side of the constrained
problem (\ref{eq:FETI-DP-system}), and to recover the primal solution.

\subsection{P-FETI-DP}

The P-FETI-DP preconditioner \cite{Fragakis-2003-MHP} is based on the
approximate solution from the first step of FETI-DP, starting from $\lambda
=0$, and with the residual $r$ as the right-hand side. The primal solution
corresponding to the result of this step is the output of the
preconditioner. Thus, from (\ref{eq:primal-sol}) with $\lambda=0$ and $f=r$,
we have%
\begin{equation}
M_{P-FETI-DP}r=E\widetilde{S}^{-1}E^{T}r.  \label{eq:operator-p-feti-dp}
\end{equation}

Comparing (\ref{eq:operator-p-feti-dp}) with the BDDC preconditioner in (\ref%
{eq:operator-bddc}), we have immediately:

\begin{theorem}
The P-FETI-DP and the BDDC preconditioners are the same.
\end{theorem}

\section{Condition Number Bounds and Eigenvalues}

\label{sec:connections}

From (\ref{eq:assembly-hat}) and (\ref{eq:operator-bddc}), the
preconditioned operator of the BDDC method is%
\begin{equation}
P_{BDDC}=\left( E\widetilde{S}^{-1}E^{T}\right) \left( R^{T}\widetilde {S}%
R\right) .  \label{eq:prec-op-bddc}
\end{equation}
From (\ref{eq:FETI-DP-system-2}) and (\ref{eq:Dirichlet-preconditioner}),
the preconditioned operator of FETI-DP\ is%
\begin{equation}
P_{FETI-DP}=\left( B_{D}\widetilde{S}B_{D}^{T}\right) \left( B\widetilde {S}%
^{-1}B^{T}\right) .  \label{eq:prec-op-feti-dp}
\end{equation}

Clearly, both the BDDC and the FETI-DP preconditioned operators have the
same general form
\begin{equation}
\left( LA^{-1}L^{T}\right) \left( T^{T}AT\right) ,
\label{eq:prec-abs}
\end{equation}%
where $A=\widetilde{S}$ is symmetric, positive definite, and
$L$ and $T$ are some linear operators such that
\begin{equation}
LT=I,  \label{eq:inv-abs}
\end{equation}%
because of (\ref{eq:BD}) and (\ref{eq:ER}). This important observation was
made in \cite{Li-2006-FBB} in the equivalent form $P\widetilde{S}^{-1}P%
\widetilde{S}:\limfunc{range}P\rightarrow \limfunc{range}$ $P$, where $P$ is
a projection, and in the present form in \cite{Brenner-2007-BFW}.

\subsection{Results for the Abstract Form of the Preconditioned Operators}

It is interesting that the fundamental eigenvalue estimate can be proved for
arbitrary operators of the form (\ref{eq:prec-abs}) -- (\ref{eq:inv-abs}).
The following lemma was proved in terms of the BDDC preconditioner in \cite[%
Theorem 25]{Mandel-2005-ATP}, and the proof carries over. Because the
translation between the two settings is time consuming, the proof (with some
simplifications but no substantial differences) is included here for
completeness. The resulting proof of the condition number bound for FETI-DP
in Theorem \ref{thm:lambda-bound} below appears to be new.

\begin{lemma}
\label{lem:lambda-bound}Let $V$ and $U$ be finite dimensional vector spaces
and $A:V\rightarrow V^{\prime}$ be an SPD\ operator. If $L:V\rightarrow U$
and $T:U\rightarrow V$ are linear operators such that $LT=I$ on $U$, then
all eigenvalues $\lambda$ of the operator $\left( LA^{-1}L^{T}\right) \left(
T^{T}AT\right) $ satisfy%
\begin{equation}
1\leq\lambda\leq\left\Vert TL\right\Vert _{A}^{2}.  \label{eq:lambda-bound}
\end{equation}
\end{lemma}

\begin{proof}
The operator $\left(  LA^{-1}L^{T}\right)  \left(  T^{T}AT\right)  $ is
selfadjoint with respect to the inner product $\left\langle T^{T}%
ATu,v\right\rangle $. So, it is sufficient to bound $\left\langle \left(
LA^{-1}L^{T}\right)  \left(  T^{T}AT\right)  u,u\right\rangle $ in terms of
$\left\langle \left(  T^{T}AT\right)  u,u\right\rangle $.

Let $u\in U$. Then
\begin{equation}
\left(  LA^{-1}L^{T}\right)  \left(  T^{T}AT\right)  u=Lw, \label{eq:def-Lw}%
\end{equation}
where $w=A^{-1}L^{T}T^{T}ATu$ satisfies%
\begin{equation}
w\in V,\quad\left\langle Aw,v\right\rangle =\left\langle T^{T}%
ATu,Lv\right\rangle \quad\forall v\in V. \label{eq:w-abs}%
\end{equation}
In particular, from (\ref{eq:w-abs}) with $v=w$ and (\ref{eq:def-Lw})%
\begin{equation}
\left\langle Aw,w\right\rangle =\left\langle T^{T}ATu,Lw\right\rangle
=\left\langle T^{T}ATu,\left(  LA^{-1}L^{T}\right)  \left(  T^{T}AT\right)
u\right\rangle . \label{eq:Aww}%
\end{equation}
and using $LT=I$, (\ref{eq:w-abs}) with $v=Tu,$ Cauchy inequality, the
definition of transpose, and (\ref{eq:Aww}),
\begin{align*}
\left\langle T^{T}ATu,u\right\rangle ^{2}  &  =\left\langle T^{T}%
ATu,LTu\right\rangle ^{2}\\
&  =\left\langle Aw,Tu\right\rangle ^{2}\\
&  \leq\left\langle Aw,w\right\rangle \left\langle ATu,Tu\right\rangle \\
&  =\left\langle Aw,w\right\rangle \left\langle T^{T}ATu,u\right\rangle \\
&  =\left\langle T^{T}ATu,\left(  LA^{-1}L^{T}\right)  \left(  T^{T}AT\right)
u\right\rangle \left\langle T^{T}ATu,u\right\rangle .
\end{align*}
Dividing by $\left\langle T^{T}ATu,u\right\rangle $, we get
\[
\left\langle T^{T}ATu,u\right\rangle \leq\left\langle T^{T}ATu,\left(
LA^{-1}L^{T}\right)  \left(  T^{T}AT\right)  u\right\rangle \quad\forall u\in
U,
\]
which gives the left inequality in (\ref{eq:lambda-bound}).

To prove the right inequality in (\ref{eq:lambda-bound}), let again $u\in U$.
Then, from (\ref{eq:def-Lw}), Cauchy inequality in the $T^{T}AT$ inner
product, definition of the $A$ norm, properties of the norm, and
(\ref{eq:Aww}),%
\begin{align*}
\lefteqn{\left\langle T^{T}ATu,\left(  LA^{-1}L^{T}\right)  \left(
T^{T}AT\right)  u\right\rangle ^{2}}\\
&  =\left\langle T^{T}ATu,Lw\right\rangle ^{2}\\
&  \leq\left\langle T^{T}ATu,u\right\rangle \left\langle T^{T}%
ATLw,Lw\right\rangle \\
&  =\left\langle T^{T}ATu,u\right\rangle \left\langle ATLw,TLw\right\rangle \\
&  =\left\langle T^{T}ATu,u\right\rangle \left\Vert TLw\right\Vert _{A}^{2}\\
&  \leq\left\langle T^{T}ATu,u\right\rangle \left\Vert TL\right\Vert _{A}%
^{2}\left\Vert w\right\Vert _{A}^{2}\\
&  \leq\left\langle T^{T}ATu,u\right\rangle \left\Vert TL\right\Vert _{A}%
^{2}\left\langle T^{T}ATu,\left(  LA^{-1}L^{T}\right)  \left(  T^{T}AT\right)
u\right\rangle .
\end{align*}
Dividing by $\left\langle T^{T}ATu,\left(  LA^{-1}L^{T}\right)  \left(
T^{T}AT\right)  u\right\rangle $, we get
\[
\left\langle T^{T}ATu,\left(  LA^{-1}L^{T}\right)  \left(  T^{T}AT\right)
u\right\rangle \leq\left\Vert TL\right\Vert _{A}^{2}\left\langle
T^{T}ATu,u\right\rangle \quad\forall u\in U.
\]

\end{proof}

The lower bound in Lemma~\ref{lem:lambda-bound} was also proved in a
different way in \cite[Lemma 3.4]{Brenner-2007-BFW}. The next abstract lemma
is the main tool in the comparison of the eigenvalues of the preconditioned
operators of BDDC and FETI-DP.

\begin{lemma}[{{\protect\cite[Lemmas 3.4 -- 3.6]{Brenner-2007-BFW}}}]
\label{Brenner} Let $V$ and $U_{i}$, $i=1,2$, be finite dimensional vector
spaces and $A:V\rightarrow V^{\prime}$ be an SPD\ operator. If $%
L_{i}:V\rightarrow U_{i}$ and $T_{i}:U_{i}\rightarrow V$ are linear
operators such that
\begin{align}
& L_{i}T_{i}=I\text{ on }U_{i},\quad i=1,2,  \label{eq:one-sided-inv} \\
& T_{1}L_{1}+T_{2}L_{2}=I\text{ on }V,  \label{eq:complementary-proj}
\end{align}
then all eigenvalues (except equal to one) of the operators $\left(
L_{1}A^{-1}L_{1}^{T}\right) \left( T_{1}^{T}AT_{1}\right) $ and $\left(
T_{2}^{T}AT_{2}\right) \left( L_{2}A^{-1}L_{2}^{T}\right) $ are the same,
and their multiplicities are identical.
\end{lemma}

\subsection{Results for FETI-DP and BDDC}

Condition number bounds now follow immediately from Lemma~\ref%
{lem:lambda-bound}.

\begin{theorem}
\label{thm:lambda-bound}The eigenvalues of the preconditioned operators of
FETI-DP and BDDC satisfy $1\leq\lambda\leq\omega_{FETI-DP}$ and $1\leq
\lambda\leq\omega_{BDDC}$, respectively, where
\begin{equation}
\omega_{BDDC}=\left\Vert E\right\Vert _{\widetilde{S}}^{2},\quad
\omega_{FETI-DP}=\Vert B_{D}^{T}Bw\Vert_{\widetilde{S}}^{2}.
\label{eq:omega}
\end{equation}
In addition, if $\widetilde{W}\neq\widehat{W}$ and (\ref{eq:corresp}) holds,
then also
\begin{equation}
\omega_{BDDC}=\omega_{FETI-DP}.  \label{eq:alt-omega}
\end{equation}
\end{theorem}

\begin{proof}
The eigenvalue bounds with (\ref{eq:omega}) follows from the form of the
preconditioned operators (\ref{eq:prec-op-bddc})\ and
(\ref{eq:prec-op-feti-dp}) and Lemma \ref{lem:lambda-bound}. The equality
(\ref{eq:alt-omega}) follows from the fact that $E$ and $B_{D}B$ are
complementary projections by (\ref{eq:corresp}), and the norm of a nontrivial
projection depends only on the angle between its range and its nullspace
\cite{Ipsen-1995-ABC}.
\end{proof}

The result in Theorem \ref{thm:lambda-bound} was proved in a different way
in \cite{Mandel-2003-CBD} for BDDC and in \cite{Mandel-2001-CDP} for
FETI-DP. For a simple proof of the bound for BDDC directly from the
variational formulation (\ref{eq:variational-BDDC}), see \cite[Theorem 2]%
{Mandel-2007-ASF}.

Equality of the eigenvalues of the two methods follows immediately from
Lemma~\ref{Brenner}:

\begin{theorem}
\label{thm:eig}Let (\ref{eq:pos-def}) -- (\ref{eq:corresp}) hold. Then, (a)
the spectra of the preconditioned operators of BDDC and FETI-DP are the same
except possibly for eigenvalue equal to one, and all eigenvalues are larger
or equal to one, and (b) the multiplicity of any common eigenvalue different
from one is the same, and the multiplicity of the eigenvalue equal to one
for FETI-DP is less than or equal to the multiplicity for BDDC.
\end{theorem}

Statement (a) of Theorem \ref{thm:eig} was proved in \cite{Mandel-2005-ATP}
in a different way, and an elegant simplified proof was given in \cite%
{Li-2006-FBB}. Statement (b) was proved in \cite{Brenner-2007-BFW}. This
presentation uses the fundamental lemma and the approach from \cite%
{Brenner-2007-BFW}.

\section{Substructuring for a Model Problem}

\label{sec:substructuring}To clarify ideas, we show how the spaces and
operators arise in the standard substructuring theory for a model problem
obtained by a discretization of the second order scalar elliptic problem.
Consider a bounded domain $\Omega\subset\mathbb{R}^{d}$ decomposed into
nonoverlapping subdomains $\Omega_{i}$, $i=1,...,N$, which form a conforming
triangulation of the domain$~\Omega$. Each subdomain $\Omega_{i}$, from now
called a substructure, is a union of Lagrangean $P1$ or $Q1$ finite
elements, and the nodes of the finite elements between the substructures
coincide. The nodes contained in the intersection of at least two
substructures are called boundary nodes. The union of all boundary nodes of
all substructures is called the interface, denoted by $\Gamma$, and\ $%
\Gamma_{i}$ is the interface of substructure$~\Omega_{i}$. The space of all
vectors of local degrees of freedom on$~\Gamma_{i}$ is denoted by $W_{i}$.
Let $S_{i}:W_{i}\rightarrow W_{i}$ be the Schur complement operator obtained
from the stiffness matrix of the substructure $\Omega_{i}$ by eliminating
all interior degrees of freedom of $\Omega_{i}$, i.e., those that do not
belong to interface $\Gamma_{i}$. We assume that the matrices $S_{i}$ are
symmetric positive semidefinite. Let $W=W_{1}\times\cdots\times W_{N}$ and
write vectors and matrices in the block form%
\begin{equation}
w=\left[
\begin{array}{c}
w_{1} \\
\vdots \\
w_{N}%
\end{array}
\right] ,\quad w\in W,\quad S=%
\begin{bmatrix}
S_{1} &  &  \\
& \ddots &  \\
&  & S_{N}%
\end{bmatrix}
.  \label{eq:block-operators}
\end{equation}
The bilinear form $a$ is then given by%
\begin{equation*}
a\left( u,v\right) =u^{T}Sv.
\end{equation*}
The solution space $\widehat{W}$ of the problem (\ref{eq:variational}) is a
subspace of $W$\ such that all subdomain vectors of degrees of freedom are
continuous across the interfaces, which here means that their values on all
the substructures sharing an interface nodes coincide.

\begin{figure}[ptb]
\begin{center}
\begin{tabular}{ccccc}
\includegraphics[width=1.55in]{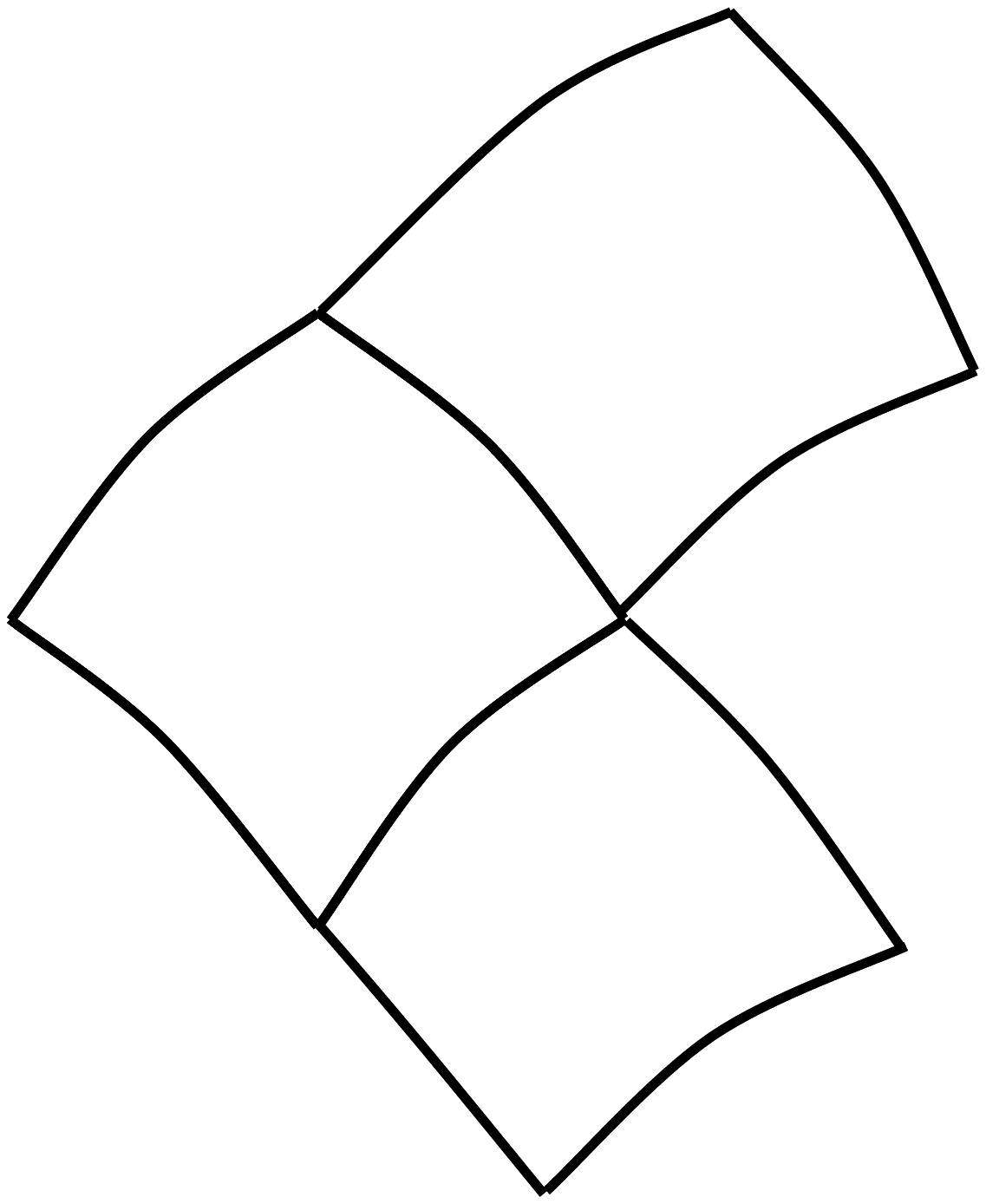} &  & %
\includegraphics[width=1.55in]{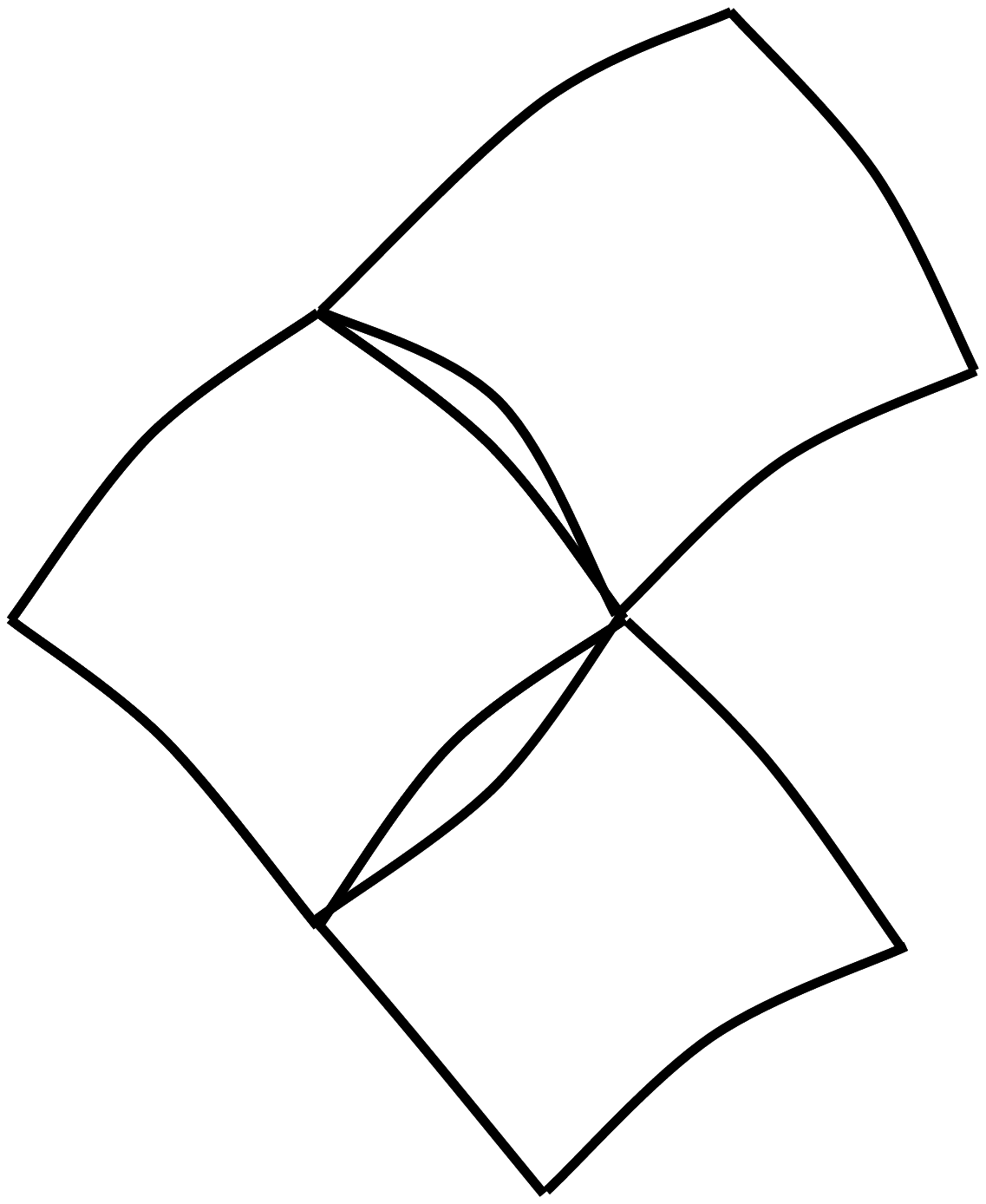} &  & %
\includegraphics[width=1.55in]{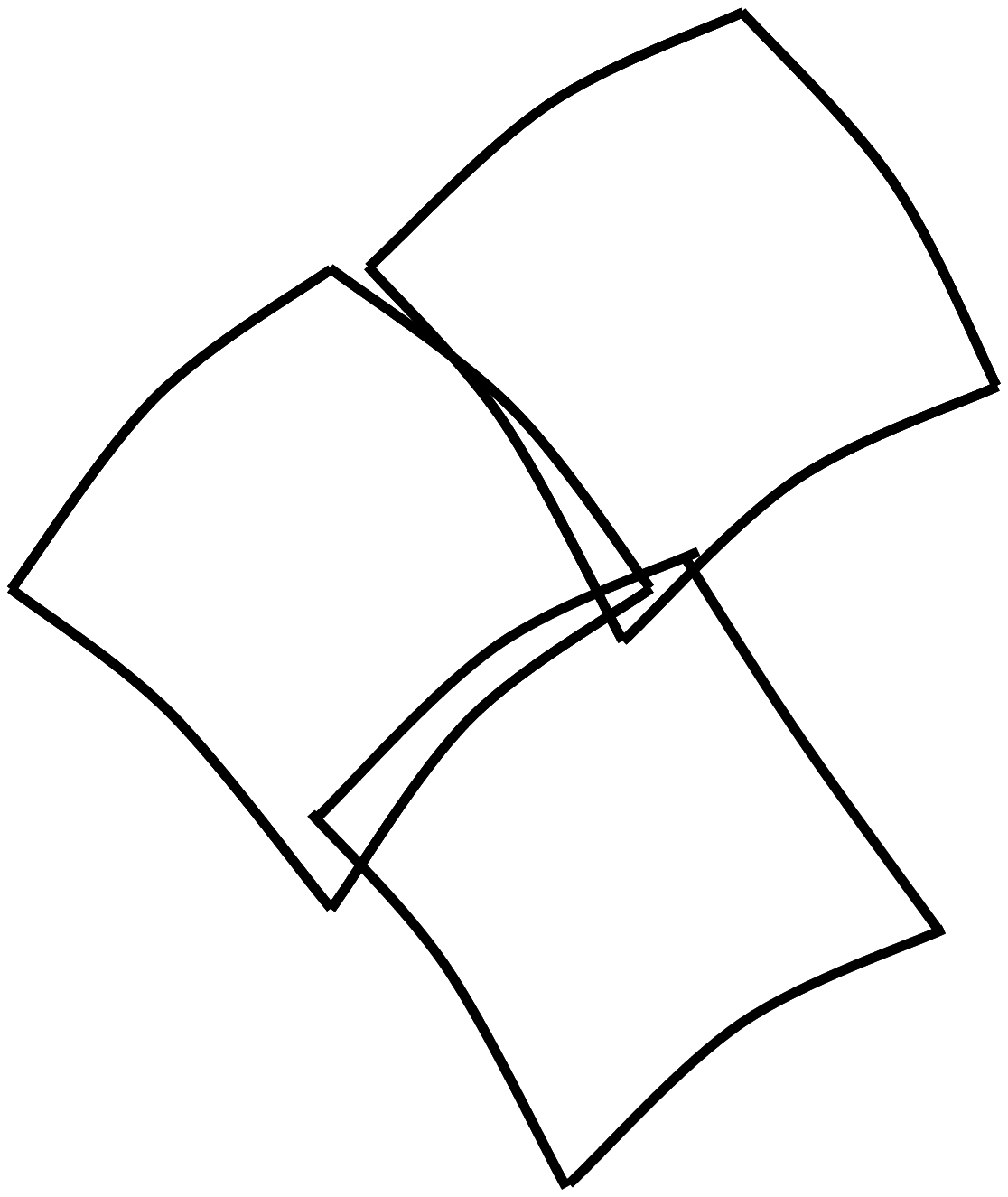} \\
$\widehat{W}$ & $\subset$ & $\widetilde{W}$ & $\subset$ & $W$%
\end{tabular}%
\end{center}
\caption{Schematic drawing of continuity conditions between substructures,
in the case of corner coarse degrees of freedom only: all degrees of freedom
continuous (the space $\widehat{W}$), only the coarse degrees of freedom
need to be continuous (the space $\widetilde{W}$), and no continuity
conditions (the space $W$).}
\label{fig:spaces}
\end{figure}

The BDDC and FETI-DP are characterized by selection of \emph{coarse degrees
of freedom}, such as values at the corners and averages over edges or faces
of subdomains (for their general definition see, e.g.,~\cite%
{Klawonn-2006-DPF}). In the present setting, this becomes the selection of
the subspace $\widetilde{W}$ $\subset W,$ defined as the subspace of all
functions such that coarse degrees of freedom are continuous across the
interfaces. Cf., Fig.~\ref{fig:spaces}. There needs to be enough coarse
degrees of freedom that the variational problem on $\widetilde{W}$ is
coercive, i.e., (\ref{eq:pos-def}) is satisfied. Creating the stiffness
matrix on the space $\widetilde{W}$ is called subassembly \cite{Li-2007-UIS}.

The last ingredients are the selections of the linear operators $E$, $B$,
and $B_{D}$. The operators $E$ and $B$ are in fact defined on the whole
space $W$; they are considered restricted on $\widetilde{W}$ only for the
purposes of the theory here. The operator $E:W\rightarrow\widehat{W}$ is an
averaging of the values of degrees of freedom between the substructure. The
averaging weights are often taken proportional to the diagonal entry of the
stiffness matrices in the substructures. The matrix $B$ enforces the
continuity across\ substructure interfaces by the condition $Bw=0$. Each row
$B$ has only two nonzero entries, one equal to $+1$ and one equal to $-1$,
corresponding to the two degrees of freedom whose value should be same. So, $%
Bw$ is the jump of the value of $w$\ between substructures. Redundant
Lagrange multipliers are possible; then $B$ does not have full row rank and $%
\Lambda =\limfunc{range}B$ is not the whole euclidean space. Finally, $B_{D}$
is a matrix such that a vector $\lambda$ of jumps between the substructures
is made into a vector of degrees of freedom $B_{D}^{T}\lambda$ that exhibits
exactly those jumps. That is, $BB_{D}^{T}=I$. The construction of $B_{D}$
involves weights, related to those in the operator $E$, so that $%
B_{D}^{T}B+E=I$, and its details are outside of the scope of this paper.
Such construction was done first for the FETI method in \cite%
{Klawonn-2002-DPF} in order to obtain estimates independent of the jump of
coefficients between substructures, and then adopted for FETI-DP. We only
note that in many cases of practical relevance, the matrix $B_{D}$ is
determined from the properties (\ref{eq:E}) -- (\ref{eq:corresp}) uniquely
as the Moore-Penrose pseudoinverse in a special inner product given by the
averaging weights in the operator $E$ \cite[Theorem 14]{Mandel-2005-ATP}.

\bibliographystyle{siam}
\bibliography{../../bibliography/bddc}

\begin{thebibliography}{10}

\bibitem{Brenner-2007-BFW}
{\sc S.~C. Brenner and L.-Y. Sung}, {\em B{DDC} and {FETI}-{DP} without
  matrices or vectors}, Comput. Methods Appl. Mech. Engrg., 196 (2007),
  pp.~1429--1435.

\bibitem{Cros-2003-PSC}
{\sc J.-M. Cros}, {\em A preconditioner for the {S}chur complement domain
  decomposition method}, in Domain Decomposition Methods in Science and
  Engineering, I.~Herrera, D.~E. Keyes, and O.~B. Widlund, eds., National
  Autonomous University of Mexico (UNAM), M\'exico, 2003, pp.~373--380.
\newblock 14th International Conference on Domain Decomposition Methods,
  Cocoyoc, Mexico, January 6--12, 2002.

\bibitem{Dohrmann-2003-PSC}
{\sc C.~R. Dohrmann}, {\em A preconditioner for substructuring based on
  constrained energy minimization}, SIAM J. Sci. Comput., 25 (2003),
  pp.~246--258.

\bibitem{Dryja-1995-SMN}
{\sc M.~Dryja and O.~B. Widlund}, {\em Schwarz methods of {N}eumann-{N}eumann
  type for three-dimensional elliptic finite element problems}, Comm. Pure
  Appl. Math., 48 (1995), pp.~121--155.

\bibitem{Farhat-2001-FDP}
{\sc C.~Farhat, M.~Lesoinne, P.~LeTallec, K.~Pierson, and D.~Rixen}, {\em
  F{ETI}-{DP}: a dual-primal unified {FETI} method. {I}. {A} faster alternative
  to the two-level {FETI} method}, Internat. J. Numer. Methods Engrg., 50
  (2001), pp.~1523--1544.

\bibitem{Farhat-2000-SDP}
{\sc C.~Farhat, M.~Lesoinne, and K.~Pierson}, {\em A scalable dual-primal
  domain decomposition method}, Numer. Linear Algebra Appl., 7 (2000),
  pp.~687--714.
\newblock {P}reconditioning techniques for large sparse matrix problems in
  industrial applications (Minneapolis, MN, 1999).

\bibitem{Farhat-1991-MFE}
{\sc C.~Farhat and F.-X. Roux}, {\em A method of finite element tearing and
  interconnecting and its parallel solution algorithm}, Internat. J. Numer.
  Methods Engrg., 32 (1991), pp.~1205--1227.

\bibitem{Fragakis-2007-FDD}
{\sc Y.~Fragakis}, {\em Force and displacement duality in {D}omain
  {D}ecomposition {M}ethods for {S}olid and {S}tructural {M}echanics}.
\newblock To appear in Comput. Methods Appl. Mech. Engrg., 2007.

\bibitem{Fragakis-2003-MHP}
{\sc Y.~Fragakis and M.~Papadrakakis}, {\em The mosaic of high performance
  domain decomposition methods for structural mechanics: {F}ormulation,
  interrelation and numerical efficiency of primal and dual methods}, Comput.
  Methods Appl. Mech. Engrg., 192 (2003), pp.~3799--3830.

\bibitem{Golub:1989:MAC}
{\sc G.~H. Golub and C.~F.~V. Loan}, {\em Matrix Computations}, Johns Hopkins
  Univ. Press, 1989.
\newblock Second Edition.

\bibitem{Ipsen-1995-ABC}
{\sc I.~C.~F. Ipsen and C.~D. Meyer}, {\em The angle between complementary
  subspaces}, Amer. Math. Monthly, 102 (1995), pp.~904--911.

\bibitem{Klawonn-2001-FNN}
{\sc A.~Klawonn and O.~B. Widlund}, {\em F{ETI} and {N}eumann-{N}eumann
  iterative substructuring methods: connections and new results}, Comm. Pure
  Appl. Math., 54 (2001), pp.~57--90.

\bibitem{Klawonn-2006-DPF}
{\sc A.~Klawonn and O.~B. Widlund}, {\em Dual-primal {FETI} methods for linear
  elasticity}, Comm. Pure Appl. Math., 59 (2006), pp.~1523--1572.

\bibitem{Klawonn-2002-DPF}
{\sc A.~Klawonn, O.~B. Widlund, and M.~Dryja}, {\em Dual-primal {FETI} methods
  for three-dimensional elliptic problems with heterogeneous coefficients},
  SIAM J. Numer. Anal., 40 (2002), pp.~159--179.

\bibitem{Li-2006-FBB}
{\sc J.~Li and O.~B. Widlund}, {\em {FETI-DP}, {BDDC}, and block {C}holesky
  methods}, Internat. J. Numer. Methods Engrg., 66 (2006), pp.~250--271.

\bibitem{Li-2007-UIS}
\leavevmode\vrule height 2pt depth -1.6pt width 23pt, {\em On the use of
  inexact subdomain solvers for {BDDC} algorithms}, Comput. Methods Appl. Mech.
  Engrg., 196 (2007), pp.~1415--1428.

\bibitem{Mandel-1993-BDD}
{\sc J.~Mandel}, {\em Balancing domain decomposition}, Comm. Numer. Methods
  Engrg., 9 (1993), pp.~233--241.

\bibitem{Mandel-2003-CBD}
{\sc J.~Mandel and C.~R. Dohrmann}, {\em Convergence of a balancing domain
  decomposition by constraints and energy minimization}, Numer. Linear Algebra
  Appl., 10 (2003), pp.~639--659.
\newblock Dedicated to the 70th birthday of Ivo Marek.

\bibitem{Mandel-2005-ATP}
{\sc J.~Mandel, C.~R. Dohrmann, and R.~Tezaur}, {\em An algebraic theory for
  primal and dual substructuring methods by constraints}, Appl. Numer. Math.,
  54 (2005), pp.~167--193.

\bibitem{Mandel-2007-ASF}
{\sc J.~Mandel and B.~Soused{\'\i}k}, {\em Adaptive selection of face coarse
  degrees of freedom in the {BDDC} and the {FETI-DP} iterative substructuring
  methods}, Comput. Methods Appl. Mech. Engrg., 196 (2007), pp.~1389--1399.

\bibitem{Mandel-2001-CDP}
{\sc J.~Mandel and R.~Tezaur}, {\em On the convergence of a dual-primal
  substructuring method}, Numerische Mathematik, 88 (2001), pp.~543--558.

\bibitem{Mangasarian-1994-NP}
{\sc O.~L. Mangasarian}, {\em Nonlinear programming}, vol.~10 of Classics in
  Applied Mathematics, Society for Industrial and Applied Mathematics (SIAM),
  Philadelphia, PA, 1994.
\newblock Corrected reprint of the 1969 original.

\bibitem{Rixen-1999-TCF}
{\sc D.~J. Rixen, C.~Farhat, R.~Tezaur, and J.~Mandel}, {\em Theoretical
  comparison of the {FETI} and algebraically partitioned {FETI} methods, and
  performance comparisons with a direct sparse solver}, Internat. J. Numer.
  Methods Engrg., 46 (1999), pp.~501--534.

\end{thebibliography}

\end{document}